\input amstex
\documentstyle{amsppt}

\NoBlackBoxes
\hcorrection{0.6in}

\define\eps{\varepsilon}
\redefine\phi{\varphi}

\define\ball{{\Bbb B^n}}
\define\ballt{{\Bbb B^2}}
\define\disk{{\Bbb D}}

\define\cdisk{{\overline \disk}}
\define\sqt{{\sqrt 2} }
\define\oosqt{{1\over \sqrt2}}
\define\dzo{{\partial \over {\partial z_1}}}
\define\dzt{{\partial \over {\partial z_2}}}
\define\dwo{{\partial \over \partial w_1}}
\define\dwt{{\partial \over \partial w_2}}
\define\hinf{{H^\infty}}
\define\hinfom{{\hinf (\Omega)}}
\define\hinfd{{\hinf (\disk)}}
\define\aom{{A (\Omega)}}
\define\normqv{{\|v\|_{H^\infty (\Omega) / I_a}}}
\define\normqvp{{\|v\|_{H^\infty (\Omega) / I_{a'}}}}
\define\normqva{{\|v\|_{A (\Omega) / J_a}}}
\define\normqfa{{\|f\|_{A (\Omega) / J_a}}}

\topmatter
\title
Finite interpolation with
minimum uniform norm in $\Bbb C^n$
\endtitle

\author
Eric Amar and Pascal J. Thomas
\endauthor
\affil
Universit\'e Bordeaux-I\\
Universit\'e Paul Sabatier
\endaffil

\abstract{Given a finite sequence $a:=\{a_1, \dots, a_N \}$ in a domain
$\Omega \subset \Bbb C^n$, and complex scalars
 $v:=\{v_1, \dots, v_N \}$,
consider the classical extremal problem
of finding the smallest uniform norm of a
holomorphic function verifying $f(a_j)=v_j$ for all $j$.
We show that the modulus of the
solutions to this problem must approach its least upper
bound along a subset of the boundary of the domain large
enough to contain the support of a measure whose hull
contains a subset of the original $a$ large enough
to force the same minimum norm. Furthermore, all the solutions
must agree on a variety which also contains this hull.
An example is given to show that the inclusions can be
strict.}
\endabstract

\thanks{Part of this work was done during a stay at the
Centre de Recerca Mathem\`atica,
Universitat Aut\`onoma de Barcelona, Bellaterra, financed
by a program of the Comunitat de Treball dels Pirineus.}
\endthanks

\endtopmatter

\document

\subheading{0. Introduction and statement of results.}

Given a finite sequence $a:=\{a_1, \dots, a_N \}$ in a domain
$\Omega \subset \Bbb C^n$, and complex scalars
 $v:=\{v_1, \dots, v_N \}$,
consider the classical extremal problem:
$$
\inf \left\{ \| f \|_\infty : f \in H^\infty (\Omega),
f(a_j)=v_j, 1 \le j \le N \right\} =: m .
\tag 1
$$
It will be convenient to consider the data $v$ as lying
in $\hinfom / I_a$, where $\hinfom$ is the algebra of
functions holomorphic and bounded on $\Omega$, and
$I_a$ is the ideal associated to $a$,
$$
I_a := \{ f \in \hinfom : f(a_j)=0, j=1, \dots,N \} .
$$
Then $m= \normqv$, by the definition
of the quotient norm.
By Montel's Theorem, we know that this problem always admits
an extremal function $f$, i.e. a representative of the class
$v$ such that $\|f\|_\infty = \normqv$.

When $n=1$ and $\Omega = \Bbb D$, it is a classical fact
that $f$ is unique, and indeed given by a constant multiple
of a Blaschke product of degree $N-1$. In particular, it is
holomorphic in a neighborhood of $\overline {\Bbb D}$ and of constant
modulus on $\partial \Bbb D$, and the same properties hold when
$\disk$ is replaced by any bounded domain $\Omega$ with smooth
boundary in the complex plane, and in even more general
one-dimensional cases (see \cite{Gr}, or
for instance \cite{Gm} and references therein).

Those properties cannot hold in higher dimension. Consider
the simple example where
$a = (\alpha,0)  \subset \disk \times \{ 0\} \subset \ball$,
the unit ball of $\Bbb C^n$,
with $\alpha:=\{\alpha_1, \dots, \alpha_N \}\subset \disk$.
Then, given any solution $f  \in H^\infty (\ball)$ to the problem (1),
the function $f_0 \in H^\infty (\disk)$ defined by
restriction ($f_0 (\zeta) := f(\zeta, 0)$) will be a solution
to the classical extremal problem in the disk, therefore
${\|v\|_{H^\infty (\ball) / I_a}}
 \ge \|v\|_{\hinfd/ I_\alpha}$.
But given any solution $f_0$
to the problem in the disk, the trivial extension $f(z):=f_0(z_1)$
will solve the problem in the ball, thus
$\|v\|_{H^\infty (\ball) / I_a}= \|v\|_{\hinfd / I_\alpha} =
\|f_0\|_\infty$,
and we have a solution given by a
constant multiple of a Blaschke product in $z_1$ which
has modulus $\|v\|_{H^\infty (\ball) / I_a}$ on $\partial \disk \times \{
0\}$ and strictly less
than $\|v\|_{H^\infty (\ball) / I_a}$ elsewhere on $\partial \ball$.

Furthermore, this solution is
not unique. Indeed, using for instance
the fact that, on the unit circle, $|f_0|^2$ can only
have critical points where $f'_0$ vanishes,
so that in this case the non-vanishing of any further derivative at that
point would imply a local violation of the maximum principle,
one sees easily that $f'_0$ does not vanish anywhere
on the circle, and thus
there exists a $\gamma>0$ such that
$|f_0 (\zeta)| \leq \|f_0\|_\infty - \gamma (1-|\zeta|^2)$
for all $\zeta \in \disk$. On the other hand, if
$(z_1,z') \in \ball \subset \Bbb C \times \Bbb C^{n-1}$,
then $|z'|^2 < 1-|z|^2$. Therefore any function of the
form $g(z_1,z') := f_0 (z_1) + \gamma h(z)$, where
$|h(z)| \le |z'|^2$, will provide another solution to
the problem.


We are
interested in the relationship between our extremal problem,
notably the sequence $a$, and
the subsets of $\partial \Omega$ where its solutions
reach their maximum modulus  $\normqv$.


\smallbreak

\proclaim{Definition}

For any function $f \in H^\infty (\Omega)$,
let
$$
M(f):= \{ \xi \in \partial \Omega :
\limsup_{z \to \xi, z \in \Omega} |f(z)| = \|f\|_\infty  \} .
$$
\endproclaim

It will be useful to highlight those subproblems of the original
problem which yield the same extremal norm.

\proclaim{Definition}
We say that $a'$
defines a {\it sufficient\/} subproblem
of (1) if and only if
$\normqvp = \normqv$.

We say that a problem $(a,v)$ is {\it minimal\/}
when it does not contain any strictly smaller sufficient
subproblem.
\endproclaim

Note that this definition depends on the values $v$.
That $a'$ be sufficient
 implies of course that the points of $a \setminus a'$
are "inactive constraints", in the sense that removing them
will not change the extremum we are looking for. Note however
that the converse is not true: it is quite possible to have
problems $(a,v)$, every constraint
$(a_i, v_i)$ of which is inactive, but of course
removing them all (resp. all but one) would lead
to a problem without constraints
whose solution is $\infty$ (resp. the modulus of the
remaining $v_j$). Take for instance three points in the disk
and values at those three points of a M\"obius
automorphism of the disk. Then any pair of points
will provide a minimal sufficient subproblem.

We denote by $A(\Omega)$ the algebra of
functions holomorphic on $\Omega$ and continuous
on $\overline \Omega$.
For any compact set $K\subset \Omega$,
the $A(\Omega)${\it-hull\/}
is defined to be
$$
 \hat K_{A(\Omega)} := \{ z \in \Bbb C^n : \forall F \in A(\Omega),
|F(z)| \leq \max_K |F| \} \ .
$$
In the case where $\overline{\Omega}$ has a
neighborhood basis of Runge domains (for instance
when $\Omega$ is convex), then we can replace $A(\Omega)$
by $\Bbb C_n [Z]$, the set of all (holomorphic) polynomials
in $n$ variables, and we just get the polynomial hull,
denoted by $\hat K$.

We will restrict attention to open sets where bounded holomorphic
functions are well approximated by functions continuous up to the
boundary, in the following sense.

\proclaim{Definition}
We say that $\Omega$ has property $(A)$ if and only if
$\Omega$ is a bounded domain and for any
$g \in \hinfom$, there exists a sequence $\{g_n\} \subset A(\Omega)$
such that for any open set $U$,
$\lim_{n \to \infty} \| g_n \|_{L^\infty (U)} = \| g \|_{L^\infty (U)}$,
and $g_n \to g$ uniformly on compacta of $\Omega$.
\endproclaim

This property holds in particular when $\Omega$ is convex and bounded
(use dilations).

\proclaim{Theorem 1}

Suppose that $\Omega$ has property $(A)$.
Let $f \in H^\infty (\Omega)$, and
$\|f\|_\infty = \normqv$. There exists a
sufficient subproblem $(a',v|_{a'})$ such that
$[M(f)]^\wedge _{A(\Omega)} \supset a'$.

In particular, if all the points of $a$ are active constraints,
then $[M(f)]^\wedge _{A(\Omega)} \supset a$, and in general
$[M(f)]^\wedge _{A(\Omega)} \cap a \neq \emptyset$.
\endproclaim

Notice that it follows from the maximum principle
that, when $v$ is not constant,
a subsequence $a'$ giving
a sufficient subproblem must contain at least two points.
In the case of the example given above,
$M(f) \supset \partial \Bbb D \times \{ 0 \}$ and
$[M(f)]^\wedge \supset \overline {\Bbb D } \times \{ 0 \}$.
In fact there is always a single set $M(f_0)$ contained
in all the $M(f)$, for any $f$ solution to the problem.

\proclaim{Lemma 2}
Given any $a$ and $v$ as above, there exists a holomorphic solution
$f_0$ to the problem (1) such that
$$
M(f_0) =
\bigcap_{f\in H^\infty(\Omega) : f(a)=v , \|f\|_\infty =\normqv}
M(f) \ .
$$
\endproclaim

In the case of the example, $M(f_0) = \partial \Bbb D \times \{ 0 \}$.
Theorem 1 says that $M(f_0)$ cannot be too small.
We give some well-known consequences in the case
when $\Omega = \Bbb B^n$.

\proclaim{Corollary 3}
\roster
\item The set $M(f_0)$ has positive (possibly infinite) length.
\item The set $M(f_0)$ cannot be
a peak-interpolation set.
\item If $M(f_0) \subset \partial \disk \times \{ 0\}$,
then $M(f_0) = \partial \disk \times \{ 0\}$.
\endroster
\endproclaim

\demo{Proof}
By applying an automorphism of the ball, we may
assume that $0 \in [M(f_0)]^\wedge \cap a$. Then \cite{Fo}
and \cite{La}
show that the length of $M(f_0)$ is at least $2 \pi$,
so remains positive after applying the inverse automorphism.

For (ii), see \cite{Ru}; (iii) is elementary.
\enddemo

\smallbreak

{\bf Representing measures}

Let $\mu$ be a Borel measure on $\overline \Omega$.
We define the {\it hull of $\mu$\/}, $\Cal E_\mu \subset \Omega$
by $z \in \Cal E_\mu$ if and only if there exists a measure
$\nu_z$, absolutely continuous with respect to $\mu$,
 which is a representing measure for $z$, i.e. for any
$f \in \aom$, $f(z)= \int_{\overline \Omega} f d \nu_z$.

\proclaim{Lemma 4}
If $\mu$ is supported on the (closed) set $K$,
that is, if $\mu (\overline \Omega \setminus K)=0$, then
$\Cal E_\mu \subset \hat K_{A(\Omega)}$.
\endproclaim

\demo{Proof}
For any $z \in \Cal E_\mu$, the measure $\nu_z$ given by
the above definition is also supported on $K$. For any
$f \in {A(\Omega)}$,
$$
|f(z)| = \left| \int_{\overline \Omega} f(\zeta) \, d \nu_z (\zeta) \right|
= \left| \int_{K} f(\zeta) \, d \nu_z (\zeta) \right|
\le \sup_K |f| \| \nu \| =  \sup_K |f|  .
$$
\enddemo

With the help of the above lemma, Theorem 1 is a consequence
of the following result.

\proclaim{Theorem 5}
Let $f_0$ be a solution to the extremal problem (1).
There exists $a'$ giving a sufficient subproblem of
$(a,v)$ and a Borel measure $\mu$
on $\overline \Omega$, supported on $M(f_0)$,
such that $a' \subset \Cal E_\mu$.
\endproclaim

The second-named author would like to thank Makhlouf Derridj,
Norman Levenberg, and Zbigniew Slodkowski for useful discussions.


\subheading{1. Proof of Theorem 5}

The methods needed to prove Theorem 5 rely
on concepts developed long ago by the first-named
author \cite{Am1}, and recently put to use to study Pick-Nevanlinna
problems in several variables \cite{Am2}. First we reduce
ourselves to $A(\Omega)$. When needed, we will write
$J_a := I_a \cap \aom$.

\proclaim{Lemma 6}
For any finite sequence $a$ and values $v$,
$\normqv = \normqva$.
\endproclaim

\demo{Proof}
Take any $f \in \hinf$ such that $f(a_j)=v_j$, $1 \le j \le N$
and $\|f\|_\infty = \normqv$.
If $\{ f_n \}$ is the sequence given by property $(A)$,
$\sup_a |f - f_n| \to 0$ as $n \to \infty$. Let $L_n$ be the
Lagrange polynomial interpolating the values $(f-f_n)(a_j)$
at the points $a_j$, then $f_n + L_n$ provide representatives
of the class of $f$ (i.e. $(f_n + L_n)-f \in I_a$) and
$\|f_n + L_n\|_\infty \to \|f\|_\infty$, so
$\normqva \le \normqv$, and the reverse inequality is trivial.
\enddemo

Now to avoid trivialities, suppose that $v \not \equiv 0$,
and let $g$ be a representative of $v$ in $A(\Omega)$.
Then, by Hahn-Banach's Theorem, there exists
$\ell \in (\aom/J_a)^*$ a continuous linear form such that
$\ell (v)= \normqva $, and
$\| \ell \| =1$; equivalently, we may consider that
$\ell \in \aom^*$, $\ell |_{I_a} \equiv 0$,
$\ell (g)= \normqva $, and $\| \ell \| =1$.

Since $\aom$ can be considered as a subspace of
$\Cal C (\partial \Omega)$, there is a measure
$\nu$ on $\partial \Omega$ which represents $\ell$,
with $\| \nu \| =1$. Write $d \nu = \theta d \mu$,
where $\mu$ is a probability measure and
$|\theta|= 1$ $\mu$-a.e.

\proclaim{Proposition 7}
Let $f_0 \in \hinf (\Omega)$ be a solution of the problem (1),
and $m= \normqv = \normqva$.
Then there exists $F^* \in L^\infty (\mu) $
(defined $\mu$-almost everywhere),
such that $| F^* | = m$, $\mu$-almost everywhere, and
$\{ \zeta \in \partial \Omega : |F^*(\zeta)|= m \}
\subset M(f_0)$.
\endproclaim

Notice that we do not prove that $F^*$ represents the
boundary values of $f_0$.

\demo{Proof}
Take a sequence $F_n := f_n + L_n$ as in the proof of
Lemma 6. By weak${}^*$ compactness of the unit ball of
$L^\infty (\partial \Omega)$, $\{ F_n \}$ admits a subsequence
which converges weakly to some $F^* \in L^\infty (\partial \Omega)$,
$\| F^* \|_\infty \le m$. Furthermore,
$$
\int F^* \theta d \mu
= \lim_{n \to \infty} \int F_n \theta d \mu
= \lim_{n \to \infty} \ell (F_n) = \ell (g) = m  \ ,
$$
so we must have $F^* (\zeta) \theta (\zeta) =m$ for
$\mu$-a.e. $\zeta$. This proves the first assertion about
$F^*$, and reduces the second one to proving that $\mu$ is supported
on the (closed) set $M(f)$.

Let $\zeta \in \partial \Omega \setminus M(f)$ and
$\psi \in \Cal C (\partial \Omega, [0,1])$ such that
$\text{supp}\, \psi \subset \overline B (\zeta,r) \cap \partial \Omega
\subset \partial \Omega \setminus M(f)$, $\psi \equiv 1$
on  $\overline B (\zeta,r/2) $. Then by definition of $M(f)$,
there exists $\eps >0$ such that
$$
\max_{\overline B (\zeta,r) \cap \partial \Omega}
|f| \le m - 2 \eps \ ;
$$
for $n$ large enough, by property $(A)$,
$\max_{\overline B (\zeta,r) \cap \partial \Omega}
|F_n| \le m - \eps$, so
$$
\left| \int \psi F_n \theta d \mu \right|
\le
(m-\eps) \int \psi d \mu ,
$$
thus
$$
\left| \int \psi F^* \theta d \mu \right|
\le
(m-\eps) \int \psi d \mu ,
$$
which implies $\int \psi d \mu =0$, thus
$\mu (\overline B (\zeta,r/2))=0$, q.e.d.
\enddemo

\smallbreak

{\bf A representation of $H^2(\mu)$}

Let $H^2(\mu)$ be the closure in $L^2 (\mu)$ of $\aom$.
Let $b \in \Omega$. Then $J_b \subset H^2 ( \mu)$; let
$e_b \in J_b^\perp \subset H^2 ( \mu)$.
Observe that $\dim J_b^\perp \le 1$.

For all
$f \in \aom$, $\langle f , e_b \rangle =
f(b) \langle \Bbb 1, e_b \rangle$. So either $J_b$
is dense in $H^2 ( \mu)$ and $e_b =0$, or
if we can find some $e_b \neq 0$,
$k_b := (\overline{ \langle \Bbb 1, e_b \rangle})^{-1} e_b$
is a reproducing kernel for the point $b$. This proves the following.

\proclaim{Lemma 8}
If $J_b^\perp \neq \{ 0 \}$, then $b \in \Cal E_\mu$.
\endproclaim

\proclaim{Definition}
For any $f \in \aom$, we let $\pi^\mu (f)$ be the
(antilinear) map from $H^2(\mu)$ to itself given by
$$
\langle \pi^\mu (f) h , k \rangle :=
\langle h, fk \rangle
$$
for any $h$, $k \in H^2 (\mu)$.
\endproclaim

\proclaim{Lemma 9} (see \cite{Am2}).

$\pi^\mu$ is an antilinear representation of $\aom$ and
$\| \pi^\mu (f) \| \le \|f\|_\infty$.

If $e_b \in J_b^\perp$, $\pi^\mu (f) (e_b) = \overline{f(b)} e_b$.
\endproclaim

Now for $s \subset \Omega$
let $E_s:= \text{span} \{ e_{b}, e_{b} \in J_{b}^\perp , b \in s\}
= J_s^\perp$, and $\pi^\mu_s$ denote the restriction of
$\pi^\mu$ to $E_s$. The above definitions could be
made for any measure $\mu$, but here we will be using the
fact that the construction of $\mu$ depended on a solution
of the problem (1).

\proclaim{Proposition 10}

For any $g \in A(\Omega)$ representing the given values $v \in
A(\Omega) / J_a$,
$\| \pi^\mu_a (g) \| = \| g \|_{\aom / J_a}
= \| g \|_{\aom / J_{a'}}$,
where $a' := \left\{ a_j \in a : J_{a_j}^\perp \ne \{0\} \right\}$.
\endproclaim

\demo{Proof}

Observe first, to avoid trivialities,
that $a'$ cannot be empty, otherwise $E_{a}$
would be reduced to $\{0\}$, and $I_{a}$
would be dense in $\hinf$, which is impossible
when $v\neq 0$,
because then we'd have solutions of the
problem (1) with arbitrarily small norm.

By the definition of $a'$, $E_a = E_{a'}$, thus
$ \pi^\mu_a   = \pi^\mu_{a'} $. Applying this to the
same function $g$, we get the same operator norms, so
it will be enough to prove the first equality to complete the proof.

Let $f$ be any function in $A(\Omega)$, $h \in J_a$.
Since the map only depends
on the values of $f$ on $a$,
$\pi^\mu_a (f+h) = \pi^\mu_a (f)$. Thus
$$
\|  \pi^\mu_a (f) \|_{op} = \|  \pi^\mu_a (f+h) \|_{op}
\le \|  \pi^\mu (f+h) \|_{op} \le \| f+h \|_\infty ,
$$
and passing to the infimum we get
$\|  \pi^\mu_a (f) \|_{op} \le \normqfa$.
So $\|  \pi^\mu_a (g) \|_{op} \le m$.

Conversely, given $g$, take $F^*$ as in Proposition 7.
Again denote $m= \normqva$. Then, since $F^*$ is obtained
as a limit of holomorphic functions, $F^* \in H^\infty (\mu)
\subset H^2 (\mu)$, and $\|F\|_{H^2 (\mu)} =m$ (because its
modulus is constant $\mu$-a.e.).

For any $h \in J_a$,
$$
\langle h, F^* \rangle = \int h \bar F^* d \mu
= \int h m \theta d \nu = m \int h d \nu = m \ell (h) = 0,
$$
by definition of $\ell$. Thus $F^* \in J_a^\perp = E_a$.

We can then test $\pi^\mu_a (g)$ on $F^*$:
$$
\langle \pi^\mu_a (g) (F^*), \Bbb 1 \rangle =
\langle F^*, g  \rangle = m \overline{\ell(g)} = m^2 .
$$
This proves the required inequality.
\enddemo

\demo{End of Proof of Theorem 5}
By Proposition 10, the subproblem defined by $a'$
is sufficient. Proposition 7 shows that the measure $\mu$ defined
after Lemma 6 is supported on $\{ \zeta: |F^*(\zeta)|=m \}
\subset M(f_0)$. And by Lemma 8, the $a'$ we have obtained
is included in $\Cal E_\mu$.
\enddemo

\subheading{2. Questions of uniqueness}

\proclaim{Definition}

The {\it uniqueness variety \/} for the problem (1) is defined by
$$
\Cal U (a,v):= \{ z \in \ball : \forall f,g \text{ solving (1) },
(f-g)(z) = 0 \} \ .
$$
\endproclaim
Clearly, $\Cal U (a,v)$ is an analytic variety containing $a$.

\proclaim{Proposition 11}
Whenever
$\mu$ is chosen as in Theorem 5, $\Cal E_\mu \subset \Cal U (a,v)$.
\endproclaim

\demo{Proof}

We reuse the notations of Lemma 6 and Proposition 7.
Suppose $f_0$ and $\tilde f_0$ are distinct solutions
to the problem (1). Take two sequences
of functions in $\aom$, $\{f_n\}$ (resp. $\{\tilde f_n\}$)
converging uniformly on compacta of $\Omega$ to $f_0$
(resp. $\tilde f_0$) and in $L^\infty (\mu)$ to $F^*$
(resp. $\tilde F^*$). The proof of Proposition 7 shows that
in fact $F^* = \tilde F^*$ $\mu$-a.e.

Suppose $b \in \Cal E_\mu$. Then, denoting by $\nu_b$ a
representing measure for $b$ that is absolutely continuous
with respect to $\mu$,
$$
f_0 (b) - \tilde f_0 (b) =
\lim_{n\to \infty} \int \left( f_n (\zeta) - \tilde f_n (\zeta) \right)
d \nu_b (\zeta) = 0
$$
by the dominated convergence theorem.
\enddemo

{\bf Examples.}
\smallbreak
In the case of the example given in the introduction,
$\Cal U(a,v) = \Bbb D \times \{ 0\} \subset \Bbb B^2$; we shall
see presently that there are some cases when $\Cal U(a,v) = \Omega$,
that is to say, the solution to the problem (1) is unique.

\proclaim{Theorem 12}

For $\Omega = \Bbb B^2$, there exists $a:=\{a_1, \dots, a_4 \}$ and $v$
such that $M(f_0)$ is a $2$-real-dimensional torus
in $\partial \ballt$, and the solution to the problem (1)
is unique.
\endproclaim

The above theorem will reduce to a result about extension
of inner functions from an analytic disk embedded into the
ball $\Bbb B^2$. First we need a simple one-variable lemma.

\proclaim{Lemma 13}

Let $a:=\{a_1, \dots, a_N \}$ be distinct points in $\Bbb D$
and $B_{N-1}$ a Blaschke product of degree exactly
equal to $N-1$. Let $v_j := B_{N-1} (a_j)$, $1 \le j \le N$.
Then $B_{N-1}$ is the unique solution to the extremal
problem (1).
\endproclaim

\demo{Proof}

It will be enough to show that $\|v\|_{\hinf(\Bbb D)/I_a}=1$. The proof
will proceed by induction. For $N=1$, $B_{N-1}$ is a
unimodular constant and the property is obvious. Suppose
it is true for $N$, and consider
$a:=\{a_1, \dots, a_N, a_{N+1} \}$.

For any $\alpha \in \Bbb D$, denote by $\varphi_\alpha$
the involutive automorphism of $\disk$ which exchanges
$0$ and $\alpha$. Suppose $f \in H^\infty (\disk)$,
$f(a)=v$, and $\|f\|_\infty <1$. Let
$g=\varphi_{B_N (a_N)} \circ f \circ \varphi_{a_N}$.
We have $g(0)=0$, so $g(\zeta)=\zeta h(\zeta)$,
with $\|h\|_\infty=\|g\|_\infty <1$.

Set $a'_j := \varphi_{a_N} (a_j)$, $v'_j := \varphi_{B_N (a_N)} (B_N (a_j))$,
$1 \le j \le N+1$. We have $v'_j = \tilde B (a'_j)$, $1 \le j \le N+1$,
where $\tilde B := \varphi_{B_N (a_N)} \circ B_N \circ \varphi_{a_N}$.
This implies that $\tilde B (\zeta) = \zeta B_{N-1} (\zeta)$,
where $B_{N-1}$ is a Blaschke product of degree $N-1$ exactly.

Now letting $v''_j := v'_j/a'_j$, $1 \le j \le N$,
we have $v''_j = B_{N-1} (a'_j)$ and
$$ \|v''\|_{\hinfd /I_{a'}}
 \le \|h\|_\infty <1 ,
$$
a contradiction with the inductive hypothesis.

\enddemo


From now on we are considering the disk embedded in the unit
ball $\Bbb B^2$ of $\Bbb C^2$ given by
$\phi (\disk) = \{ \phi (\zeta): \zeta \in \disk \}$ where
$\phi (\zeta) := \oosqt (\zeta , \zeta^2)$.
Observe that
$\phi (\disk) = \{ (z_1,z_2) \in \ball : z_2 = \sqt z_1^2  \}$.

\proclaim{Lemma 14}
Suppose that $g \in \Cal H (\cdisk)$ is an inner function
(i.e. $|g(e^{i\theta})|=1$ for all $\theta \in \Bbb R$),
analytic in a neighborhood of the closed unit disk,
such that there exists $\tilde g \in H^\infty (\ballt)$
with $\tilde g (\phi (\zeta))= g(\zeta)$ for all $\zeta \in \disk$
and $\| \tilde g \|_\infty = \| g \|_\infty =1$.
Then $g'(0)=0$, and if we write $g(\zeta) = \zeta h (\zeta)$,
there exists $H$ holomorphic in $\ballt$ such that
$$
\tilde g (z_1 , z_2)= g( \sqt z_1) +
(z_2 - \sqt z_1^2) {\sqt \over 3} h(  \sqt z_1 )
+ (z_2 - \sqt z_1^2)^2 H(z_1 , z_2) \ .
$$

\endproclaim

\demo{Proof}

Step 1:
Claim

For any differentiable function $f$ on the ball, set
$$
Lf(z) := \left( \dzo f - \sqt z_1 \dzt f \right) (z) \ .
$$
Then for any $\zeta \in \disk$, $L \tilde g (\phi(\zeta)) =0$.

This is to be compared with \cite{Ru,Theorem 11.4.7}.

\smallbreak

If $\tilde g$ was assumed to be smooth in a neighborhood
of $\overline \ballt$, it would be enough to notice that
for each $z \in \phi (\partial \disk) \subset \partial \ballt$,
$L$ is a derivation along the complex tangent line
to $\partial \ballt$ at $z$. Since $|\tilde g|$ is maximal
on $\phi (\partial \disk)$ with respect to $\partial \ballt$,
its derivative $L \tilde g$ should vanish there. The slightly
more intricate argument that follows merely extends this
to the case where $\| \tilde g \|_\infty  = 1$.

Notice first that since $g$ is smooth across the unit circle,
its derivative is bounded in a neighborhood of it and we have
$c_1>0$ such that $1- |g(\zeta)|^2 \le c_1(1-|\zeta|^2)$ for
all $\zeta \in \disk$. Now consider the complex line $\Cal L$
passing through the point $\phi (\zeta_0)$ and parallel
to the vector $(1,-\zeta_0)$. Since when $|\zeta_0|=1$
this is the complex tangent to $\partial \ballt$ at $\phi (\zeta_0)$,
there exists a $c_2>0$ such that the disk of
center $\phi (\zeta_0)$, of radius
$c_2 (1-|\zeta_0|^2)^{1/2}$ along the line $\Cal L$ is contained
in $\ballt$. Then the function
$$
f(\zeta):= \tilde g \left( \phi (\zeta_0) +
c_2 (1-|\zeta_0|^2)^{1/2} \zeta (1,-\zeta_0) \right)
$$
is bounded by $1$ in modulus on the unit disk, and Schwarz-Pick's
Lemma (see \cite{Gr, Chap. I, Lemma 1.2}) shows that
$$
|f'(0)| \le 1-|f(0)|^2 = 1 - |\tilde g (\phi(\zeta_0))|^2
= 1- |g(\zeta_0)|^2 \le c_1(1-|\zeta_0|^2) \ ,
$$
and since $f'(0)= c_2 (1-|\zeta_0|^2)^{1/2} L\tilde g (\phi(\zeta_0))$
(notice that along $\phi (\disk)$, $z_1 \sqt = \zeta$),
we have $|L\tilde g (\phi(\zeta_0))| \le C (1-|\zeta_0|^2)^{1/2}$.
Now $L\tilde g (\phi(\zeta))$ is a holomorphic function on $\disk$,
so it must be identically zero, which proves the Claim.

Step 2.

Consider the change of variables
$$
\left\{
\aligned w_1 &= z_1 \\
w_2 &= z_2 - \sqt z_1^2 \endaligned
\right\}
\quad
\Leftrightarrow
\quad
\left\{
\aligned z_1 &= w_1 \\
z_2 &= w_2 + \sqt w_1^2 \endaligned
\right\}
\ .
$$
If we set $\tilde g_1 (w_1,w_2):= \tilde g (z_1,z_2)$, we then
have $L_1 \tilde g_1 (w_1,0) =0$, where
$L_1 \tilde g_1 (w) =
\left( \dwo \tilde g_1 - 3 \sqt w_1 \dwt \tilde g_1 \right) (w)$.

Since $\tilde g_1 (w_1,0) = g (\sqt w_1 )$, we have
$\dwo \tilde g_1 (w_1,0) = \sqt g' (\sqt w_1 )$, and the above
partial differential equation becomes
$$
\sqt g' (\zeta ) = 3 \zeta \dwt \tilde g_1 ({\zeta \over \sqt},0) ,
$$
which can be solved if and only if $g'(\zeta) = \zeta h(\zeta)$.
We then have $\dwt \tilde g_1 (w_1,0) = {\sqt \over 3} h(\sqt w_1)$,
and this provides the expansion of order $1$ of $\tilde g_1$
near $\{ w_2 =0 \}$, so
$$
\tilde g_1 (w_1,w_2) = g (\sqt w_1)+ w_2 {\sqt \over 3} h(\sqt w_1)
+ w_2^2 H_1 (w_1,w_2) \ ,
$$
for some $H$ holomorphic on the image of $\ballt$ under the
change of variables. Going back to the $(z_1,z_2)$ variables,
we get the Lemma.
\enddemo
\smallbreak

We now make a small aside to look into the problem of
extending a family of simple functions (the monomials $\zeta^k$)
from the analytic disk $\phi (\disk)$ to $\ballt$ with the smallest
possible $H^\infty$ norm.

\proclaim{Lemma 15}

(i) For any $\tilde g \in \Cal H (\ballt)$ such that
$\tilde g (\phi(\zeta))=\zeta$, $\| \tilde g \|_\infty \ge \sqt$,
and this bound is attained by $\tilde g (z_1, z_2) = \sqt z_1$.

(ii) For any $k \ne 1$, there exists $\tilde g_k \in \Cal H (\ballt)$
such that $\tilde g_k (\phi(\zeta))=\zeta^k$ and
$\| \tilde g_k \|_\infty =1$.

In particular, one can take
$\tilde g_2 (z_1,z_2)= {2\over 3} (z_1^2 + \sqt z_2)$,
$\tilde g_3 (z_1,z_2)= 2 z_1 z_2$.

\endproclaim

{\bf Remark.}

We know that the only analytic
disks in the ball that allow the uniform norm-preserving
extension of any bounded holomorphic function are the affine
embeddings of $\disk$ \cite{St}, \cite{Su}; this is
to be compared with Lempert's result that the only disks which admit
a holomorphic retraction are geodesic disks for the Kobayashi distance,
i.e. in this instance affine disks once again
\cite{Le1}, \cite{Le2}.  For our disk
$\phi (\disk)$, the above Lemma gives explicit examples of the
functions which do or don't admit norm-preserving extensions.

\demo{Proof}

(i) Since $\tilde g({\zeta \over \sqt} , {\zeta^2 \over \sqt}) = \zeta$,
${1 \over \sqt} \dzo \tilde g (0,0) = 1$. Applying Schwarz's Lemma, we get
$\sup_{z \in \disk} \tilde g ( z,0) \ge \sqt$, whence the result.

(ii) When $g(\zeta) = \zeta^2$, $g'(\zeta) = 2 \zeta$,
$h(\zeta) = 2$, and setting $H=0$, we find $\tilde g_2$. Checking
the norm inequality for $(z_1,z_2) \in \partial \ballt$ is
elementary.

In the same way, or by inspection, we find $\tilde g_3$. Given
any integer $k \ge 2$, we can find two
non negative integers $a$, $b$ such that
$k = 2a + 3b$. We then set $\tilde g_k = \tilde g_2^a \tilde g_3^b$.

\enddemo

The next result will essentially complete the proof of Theorem 12.

\proclaim{Lemma 16}

The function $\tilde g_3$ in Lemma 14 is the only
 $\tilde g \in H^\infty ( \ballt)$ such that
$\tilde g (\phi(\zeta))=\zeta^3$ and $\| \tilde g \|_\infty =1$.

\endproclaim

\demo{Proof}

Let us first consider the simpler case where
$\tilde g$ is holomorphic in a neighborhood of the closed ball.
Then so is $H$ (the function obtained in Lemma 14).

For any $\nu \in ]0,2\pi[$, consider the map from the disk to the
ball given by
$\psi_\nu (\zeta) := ({\zeta \over \sqt}, e^{i\nu} {\zeta^2 \over \sqt} )$.
Then, applying Lemma 14,
$$
\tilde g (\psi_\nu (e^{i\theta}) ) =
e^{i\nu} e^{3 i\theta} + {1 \over 2 }  e^{4 i\theta} (e^{i\nu} -1)^2
H (\psi_\nu (e^{i\theta}) ) \ .
$$
A winding number argument then shows that  $H$ must vanish at
some point
along the curve $\psi_\nu ( \partial \disk )$. We will show
that in fact $H$ is identically zero.

Now remove the additional assumption;
$$
\tilde g (\psi_\nu (\zeta) ) =
e^{i\nu}  \zeta^3  \left[ 1 + {1 \over 2 }  \zeta (e^{i\nu} -1)^2
H (\psi_\nu (\zeta) ) \right] \ .
$$
Set $f(\zeta) = {1 \over 2 }  \zeta (e^{i\nu} -1)^2
H (\psi_\nu (\zeta) ) $. This function can only be
constant if $H = 0$. Suppose this is not the case.

We claim that for any $r\in ]0,1[$,
there exists $\theta_r \in ]0,2\pi[$ such that
$f(r e^{i\theta_r}) >0$. Indeed, there exists $\delta > 0$ such that
this is true for all $r\in ]0,\delta[$, by the Open Mapping Theorem,
since $f(0)=0$. Let $r_0$ be the largest number such that
the conclusion of the claim holds for all $r \in ]0,r_0[$.
If $r_0 <1$, since the winding number of the curve
$f(re^{i \theta})$ around $0$ is positive for $r$ small enough
and can only change for a value of $r$ at which the curve goes through
$0$, there must be $0<r_1\le r_0$ such that
$f(r_1 e^{i\theta_{r_1}}) =0$. Then there is $r_2 \in ]0,r_1[$
such that $f(r_2 e^{i\theta_{r_2}})$ is maximal (we use the
compactness of $f[\overline D (0,r_1)]$). But this violates the
Open Mapping Theorem in a
neighborhood of the point $r_2 e^{i\theta_{r_2}}$.

By the same argument, we can see that the function which to
$r$ associates the largest possible value $f(r e^{i\theta_r}) $
cannot have a local maximum, and that (with a slight
abuse of notation)
$\limsup_{r\to 1^-} f(r e^{i\theta_r}) > 0$.
This yields
$$
\limsup_{r\to 1^-}  |\tilde g (\psi_\nu (r e^{i\theta_r}) ) | \ge
\limsup_{r\to 1^-} r^3  \left[ 1 + f(r e^{i\theta_r}) \right] >1 ,
$$
a contradiction. Thus we must have $H \equiv 0$, q.e.d.
\enddemo

\smallbreak

\demo{End of Proof of Theorem 12}

Pick $a_j := \phi (\zeta_j)$ where $\zeta_1 , \dots , \zeta_4$
are distinct points in the unit disk, and $v_j := \zeta_j^3$.
By Lemma 12, any solution to the problem (1) must be equal
to $\zeta^3$ at the point $\phi (\zeta)$, for any $\zeta \in \disk$,
and $\| v \|_{\hinf (\ballt)/I_a} \ge 1$. By Lemma 15,
$\| v \|_{\hinf (\ballt)/I_a} = 1$, and by Lemma 16,
the solution to the problem is unique and assumes its maximum
modulus on the set $\{ |z_1|=|z_2|= {1 \over \sqt} \}$.
\enddemo

This example shows that the inclusions proved in
Theorem 1 and Proposition 11 can be strict. Here
$f_0(z) = 2 z_1 z_2$, $\Cal U (a,v) = \ballt$,
and
$M(f_0)= \{ |z_1|=|z_2|= {1 \over \sqt} \}$,
so $ M(f_0)^\wedge=\{ |z_1| \le {1 \over \sqt},
|z_2|\le {1 \over \sqt} \}$.

On the other hand, since $f_0 \in A(\ballt)$
already, we have $F^*=f_0 |_{\partial \ballt}$.
It is elementary to see that
the form $\ell$ can be represented by integration
against a function along the boundary of the
embedded disk $\varphi (\disk)$, so we may take
$\mu = \varphi_* (\frac1{2\pi} d \theta)$,
and $\Cal E_\mu = \varphi (\disk)$. So we have
the strict inclusions that we had announced.

\vfill\eject

\bigskip

{\it Keywords}: analytic discs, extremal problems,
extension of analytic functions,
Pick-Nevanlinna.

\medskip

1991 AMS Subject Classification:

Primary 32E30

Secondary 32A35, 30D50, 30C80, 32E20.

\bigskip

Eric Amar

U.F.R. Math\'ematiques

Universit\'e Bordeaux-I

351 cours de la Lib\'eration

33405 TALENCE, France

 e-mail: eamar\@math.u-bordeaux.fr

\smallskip

Pascal J. Thomas

Laboratoire Emile Picard

Universit\'e Paul Sabatier

 118 route de Narbonne

31062 Toulouse
Cedex, France

e-mail: pthomas\@cict.fr

\end